\documentclass{article}

\usepackage{arxiv}

\usepackage{tikz}
 \usetikzlibrary{calc,trees,positioning,arrows,chains,shapes.geometric,%
decorations.pathreplacing,decorations.pathmorphing,shapes,%
matrix,shapes.symbols,plotmarks,decorations.markings,shadows}

 \definecolor{mygray}{gray}{0.8}
 \definecolor{mmygray}{gray}{0.6}    
    
\usepackage{tikz}
\usepackage{verbatim}
\usepackage{subfig}

\usepackage[T1]{fontenc}    % use 8-bit T1 fonts
\usepackage{hyperref}       % hyperlinks
\usepackage{url}            % simple URL typesetting
\usepackage{booktabs}       % professional-quality tables
\usepackage{amsfonts}       % blackboard math symbols
\usepackage{nicefrac}       % compact symbols for 1/2, etc.
\usepackage{microtype}      % microtypography
\usepackage{amsmath}
\usepackage{amssymb}
\usepackage{amsmath, amsthm}
\usepackage{graphicx}
\usepackage[colorinlistoftodos]{todonotes}
\usepackage{datetime}
\theoremstyle{definition}

\title{Extension of Zorn's lemma to arbitrary binary relations}

\author{
  Athanasios Andrikopoulos\thanks{Professor  (https://www.ceid.upatras.gr/webpages/faculty/aandriko/)} \\
  Dept. of Computer Engineering and Informatics\\
  University of Patras\\
  Patras, 26504, Greece \\
  \texttt{aandriko@ceid.upatras.gr} \\
}

\begin{document}
\maketitle

In 1935, Max Zorn introduced Zorn's lemma with the intention of shortening proofs 
in algebra that had previously used the Axiom of Choice or the Well-Ordering theorem.
However, there are several applications in optimization, economics, decision analysis, 
and game theory where a binary relation is not transitive as required by Zorn's lemma.
For example,
a choice process in social choice theory is formally modeled as a maximization problem. 
More specifically, a choice process is presented by a choice function that associates with 
each feasible set of alternatives a certain subset of it, which consists of the maximal 
elements according to the viewpoint of a binary relation. However, the set of maximal 
elements is often empty. In this case, the crucial question which arises is what to count as a choice. 
That is, what sets of alternatives may be considered as reasonable solutions? To answer this question, 
a number of theories, called general solution theories, have been
proposed to take over the role of maximality in the absence of maximal elements.
Zorn's lemma cannot be applied to solve these problems 
because binary relations are not transitive.
In this note, the Zorn lemma is extended to arbitrary binary relations and thus
the Zorn lemma can do for optimization when the transitivity is broken. 
Zorn's extended lemma can be used to prove existence theorems of generalized solution concepts  
for binary relations that do not have maximal elements.

\par\smallskip\noindent
{\bf Notations and definitions}
Let $X$ be a (finite or infinite) non-empty set of alternatives, and let $R \subseteq X \times X$ be a binary relation on $X$. 
We say that $R$ on $X$ is (i) {\it reflexive} if for each $x \in X$, $(x, x) \in R$; (ii) {\it transitive} if for all 
$x, y, z \in X,[(x, z) \in R$ and $(z, y) \in R] \Rightarrow(x, y) \in R$; (iii) {\it antisymmetric} if for each $x, y \in X$, $[(x, y) \in R$ 
and $(y, x) \in R] \Rightarrow x=y$. A binary relation $\preceq$ on $X$ is a {\it partial order} if it is reflexive, 
transitive and antisymmetric. 
 A {\it total order} $\preceq$ on $X$ is a partial order in which any two elements are comparable.
A well-order $\preceq$ on a set $X$ is a total order on $X$ with the property that every non-empty subset of $X$
has a least element with respect to $\preceq$. 
Given a binary relation $R$, the asymmetric part $P(R)$ of $R$ is defined as follows:
\begin{center}
$P(R)=\{(x, y) \in X \times X \mid(x, y) \in R \text { and }(y, x) \notin R\} .
$
\end{center}
The {\it transitive closure} of $R$ is denoted by $\overline{R}$, that is for all $x, y \in X,(x, y) \in \overline{R}$ if 
there exist $K \in \mathbb{N}$ and $x_{_0}, \ldots, x_{_K} \in X$ such that 
$x=$ $x_{_0},\left(x_{_{k-1}}, x_{_k}\right) \in R$ for all $k \in\{1, \ldots, K\}$ and $x_{_K}=y$. 
If $\approx$ is an equivalence relation on $X$, then the {\it quotient set} by this equivalence relation $\approx$
will be denoted by $X_{_\approx}$,
and its elements (equivalence classes) by $[x]$.
If $A\subseteq X$, then $[A]=\{[x]\vert x\in A\}$.
The surjective map $\pi: x\longrightarrow [x]$ from $X$ onto $X_{_\approx}$, 
which maps each element to its equivalence class, is called the {\it canonical projection}. 
Let $[x]\in \pi(A)\subseteq X_{_\approx}$, then
since $\pi$ is a surjection map
from $X$ onto $X_{_\approx}$, for each $[x]\in [A]$ there exists an element $x\in A$ such that $x=\pi^{-1}([x])$. 
When such an element is chosen, it is called a {\it representative} of the class. 
We define an equivalence relation $\approx$ on $X$
\begin{center}
$x\approx y$ if and only if $x=y$ or $(x,y)\in \overline{R}$ and $(y,x)\in \overline{R}$.
\end{center}
Because $(x,y)\in \overline{R}$ is an equivalence relation, it is easy to verify that $R$ defines a partial order
$R^{\ast}$ on $X_{_\approx}$ by:
\begin{center}
$[x]R^{\ast}[y]$ if and only if $[x]=[y]$ or $[x]\neq [y]$ and there exists $x^{\prime}\in [x]$ and $y^{\prime}\in [y]$ such that $x^{\prime}\overline{R} y^{\prime}$.
\end{center}
The relation $R^{\ast}$ is called the {\it quotient relation of} $R$ in space $X_{_\approx}$. Clearly, $R^{\ast}$ is a partial order.
If $R$ is a binary relation on $X$, then $(X_{_\approx}, R^{\ast})$ is called the {\it quotient space} of
the space $X$ endowed with the binary relation $R$.

An $R$-{\it chain}, denoted $\mathfrak{C}$, is a subset of $X$ such that 
$x, x^{\prime} \in \mathfrak{C}$ implies $xR x^{\prime}$ 
or $x^{\prime}R x$. 
Note that the empty set is an $R$-chain.
An element $x\in X$ is an $R$-{\it upper bound} (resp. a {\it strict} $R$-{\it upper bound})
of a subset $B$ of $X$
if $xRy$ (resp. $(xP(R)y$) for all $y\in B$.
An element $x\in X$ is said to be $R$-{\it maximal} if for all $y\in X$ it holds that $(y,x)\notin P(R)$. 
A subset $Y \in X$ is $R$-{\it undominated} if and only if for no $x \in Y$ there is a $y \in X \backslash Y$ such that $yRx$.
If $Y=\{x\}$, then $\{x\}$ is an $R$-undominated element.
An $R$-undominated set is a {\it minimal $R$-undominated} if none of its proper subsets has this property.
A subset $Y \subseteq X$ is an $R$-{\it cycle} if, for all $x, y \in Y$, we have $(x, y) \in \overline{R}$ and $(y, x) \in \overline{R}$. 
The $R$-cycle $Y$ is a {\it top} $R$-{\it cycle} if there is no $y \in X \backslash Y$ and there is no $x \in Y$ 
such that $(x, y) \in R$. 
The {\it strong top} $R$-{\it cycle} is a top $P(R)$-cycle. 
Any minimal $R$-undominated set is an $R$-undominated element or a strong top $R$-cycle (see \cite{and}).
\par\smallskip\noindent
{\bf Axiom of Choice}: A {\it choice function} on a set $X$ is a function $f:2^{^X}\setminus \emptyset\longrightarrow X$
such that $f(A)\in A$ for every non-empty $A\subset X$. The {\it Axiom of Choice} asserts that on every set there is a choice function.

In partially ordered
sets the notions of $\preceq$-undominated element and $\preceq$-maximal element coincide.
If $C$ is an $\preceq$-chain in $X$ and $x\in C$, then we define (see \cite{lew})
\begin{center}
$Q(C,x)=\{y\in C\vert x>y\}$.
\end{center}

Given a partial order $\preceq$ in $X$, by using the axiom of choice (see \cite{lew}), 
we choose a function $f$ that assigns to every $\preceq$-chain 
$\mathfrak{C}$ of $X$, a strict upper bound $f(\mathfrak{C})$.
We shall say that a subset $A$ of $X$ is {\it conforming} if the following two conditions hold:

(a) $A$ is a $\preceq$-chain such that every subset of $A$ has a least element with respect to $\preceq$;

(b) For every element $x\in A$, we have $x=f(Q(A,x))$.

A social choice function is a functional relationship, 
$\mathcal{C}: \mathcal{P}(X) \rightarrow \mathcal{P}(X)$ such that, for every 
$A \in \mathcal{P}(X), \mathcal{C}(A)$ is a non-empty subset of $A$ which 
represents those alternatives chosen by the individual or society.

\par\smallskip\noindent
{\bf Schwartz set}:
{\it Generalized Optimal-Choice Axiom} ($(\mathcal{G}\mathcal{O}\mathcal{C}\mathcal{H} \mathcal{A}$)
(Schwartz (1986), [Page 142]). For each $A\subseteq X$, $\mathcal{C}(A)$
is equivalent to the union of all minimal $R$-undominated subsets of $A$. 
The {\it Schwartz set} is the choice set from a given set specified by the 
$(\mathcal{G}\mathcal{O}\mathcal{C}\mathcal{H} \mathcal{A}$)
condition.That is, for each $A\subseteq X$, $\mathcal{C}(A)=\displaystyle\bigcup_{_{B\in \mathcal{D}}}B$, 
where $\mathcal{D}$ is the set of all minimal $R$-undominated subsets
of $A$. Equivalently, the Schwartz set is the union of all $R$-maximal elements and all strong top $R$-cycles in $X$.

\par\smallskip\noindent
{\bf Zorn's lemma}: If every chain of a partially ordered set $X$ has an upper bound, then $X$ 
has a maximal element.

\par\smallskip\noindent
 The following lemma is proved in \cite{rud}.
\par\smallskip\noindent
{\bf Lemma}. Suppose $\mathcal{F}$ is a nonempty collection of subsets of a 
set $X$ such that the union of every subchain of $\mathcal{F}$ belongs to $\mathcal{F}$. 
Suppose $g$ is a function which associates to each $D \in \mathcal{F}$ a set $g(D) \in \mathcal{F}$ 
such that $D \subset g(D)$ and $g(D)-D$ consists of at most one element. 
Then, there exists a $D \in \mathcal{F}$ for which $g(D)=D$.

\par\smallskip\noindent
{\bf Main result} We now present the main result of this note, which contains  the following theorem. 

\par\bigskip\noindent
{\bf Theorem}.
If $R$ is a binary relation on a set $X$, then the following statements are equivalent:
\par
(i) For every set $X$ there is a choice function.
\par
(ii) If each $R$-chain in $X$ has an $R$-upper bound, 
then $X$ has a top $R$-cycle.
\par\smallskip\noindent
{\it Proof}. (i)$\Rightarrow$(ii).\ 
Let $\mathfrak{C}$ be the collection of all  $R$-chains of $X$. 
Since an arbitrary element of $X$ serves as 
an upper bound for the empty $R$-chain, the family $\mathfrak{C}$ is not empty. 
Let $f$ be a choice function for $X$. If $D \in \mathfrak{C}$, let $D^\ast$ be the set of all $x$ 
in the complement of $D$ such that 
$x$ is an $R$-upper bound of $D$.

We first prove that $D$ contains a minimal $R$-undominated subset $\widetilde{D}$.
If $D=\emptyset$, then $D$ is itself an $R$-minimal undominated set.
Suppose that $D\neq\emptyset$.

 Put
\begin{center}
$g(D)=D \cup\{f(D^\ast)\}$,\ if $D^\ast \neq \emptyset$  
and
$g(D)=D$,\ if $D^\ast=\emptyset$.
\end{center}
Then, $g(D)=D\cup \{x^{\ast}\}$, for some $x^{\ast}\in D^{\ast}$, belongs to $\mathfrak{C}$. On the other hand,
$D \subset g(D)$ and $g(D)-D$ consists of at most one element.
By the lemma above, $f(D^\ast)=\varnothing$ or $D^\ast=\varnothing$ for at least one $D \in \mathfrak{C}$.  
It  follows that $D$ is an $R$-undominated subset of $X$.
We prove that $X$ has a minimal $R$-undominated subset. If this is $D$, we have nothing to prove.
Otherwise, there exists at least a $t_{_1}\in D$ such that $D\setminus \{t_{_1}\}$ is an $R$-undominated set in $X$.
Hence, $(t_{_1},x)\notin R$ for each $x\in D\setminus \{t_{_1}\}$. It follows that $(t_{_1},x)\notin \overline{R}$ for each $x\in D\setminus \{t_{_1}\}$.
Indeed, suppose to the contrary that $(t_{_1},x)\in \overline{R}$ for each $x\in D\setminus \{t_{_1}\}$.
It then follows that, there exists a natural number $n$ and alternatives $\mu_{_1},\mu_{_2},...,\mu_{_n}$ such that
$t_{_1}R\mu_{_1}R\mu_{_2}R...R\mu_{_n}Rx$.
Therefore, $\mu_{_n}\in D\setminus \{t_{_1}\}$, for suppose otherwise: since $x\in D\setminus \{t_{_1}\}$, we cannot have $\mu_{_n}Rx$.
Similarly, $\mu_{_{n-1}}\in D\setminus \{t_{_1}\}$, and an induction argument based on this logic yields $t_{_1}\in D\setminus \{t_{_2}\}$,
a contradiction.
Since $D$ is an $R$-chain, we conclude that $(x,t_{_1})\in R\subseteq \overline{R}$ and thus $(x,t_{_1}) P(\overline{R})$.
Similarly, if $D\setminus \{t_{_1}\}$ is not a minimal $R$-undominated subset of $X$, there exists $t_{_2}\in D\setminus \{t_{_1}\}$ 
such that $D\setminus \{t_{_1},t_{_2}\}$ is a $R$-undominated subset of $X$ and  $(t_{_2},t_{_1}) P(\overline{R})$ and so on.
Let $\Gamma=(t_i)_{_{i\in I}}$ be a net in $D$ such that $t_{_j}P(\overline{R})t_{_i}$ for each $j>i$ and $D\setminus \Gamma$ 
is an $R$-undominated set in $X$. Clearly, $\Gamma$ is an $P(\overline{R})$-chain 
in $D$. 
Let $\mathcal{G}$ be the set of all $P(\overline{R})$-chains $G$ in $D$ such that
$D\setminus G$ 
is an $R$-undominated set in $X$. 
If $G\in \mathcal{G}$, let $G^{\ast}$
be the set of all $x$ in $D\setminus G$ such that $G\cup \{x\}\in \mathcal{G}$. Since $\{t_{_1}\}\in \mathcal{G}$, 
this family is non-empty.
Let $\mathfrak{f}$ be a choice function for $D$.
 Put
\begin{center}
$h(G)=G \cup\{\mathfrak{f}(G^\ast)\}$,\ if $G^\ast \neq \emptyset$  
and
$h(G)=G$,\ if $G^\ast=\emptyset$.
\end{center}
Then, $h(G)=G\cup \{x^{\ast}\}$, for some $x^{\ast}\in G^{\ast}$, belongs to $\mathcal{G}$. On the other hand,
$G \subset h(G)$ and $h(G)-G$ consists of at most one element.
By the lemma above, $G^{\ast}=\emptyset$ 
for at least one $\widetilde{G}\in \mathcal{G}$, and any such $\widetilde{G}$ is a
maximal element of $\mathcal{G}$.
It follows that $\widetilde{D}=D\setminus \widetilde{G}$ is a minimal $R$-undominated subset of $X$.
Indeed,
Suppose to the contrary that 
there exists at least one $\lambda\in \widetilde{D}$ such that
$\widetilde{D}\setminus \{\lambda\}$ is an $R$-undominated subset of $X$. 
Let $E=\widetilde{G}\cup \{\lambda\}$.
Then, 
since $(t,\lambda)\in P(\overline{R})$ we have that $(t,x)\in P(\overline{R})$ for each $x\in\widetilde{D}\setminus E$ and
$D\setminus E$ 
is an $R$-undominated subset of $X$. Hence, $E\in \mathcal{G}$, a contradiction to the maximal character of $\widetilde{G}$.
The last contradiction shows that $\widetilde{D}$ is a minimal $R$-undominated subset of $X$. 
We have two cases to consider: (i) $\widetilde{D}\neq \emptyset$ and (ii)  $\widetilde{D}=\emptyset$.
\par\smallskip\par\noindent
{\it Case {\rm (i)}}: We prove that 
$\widetilde{D}$ is a top $R$-cycle.
Suppose that $x \in \widetilde{D}$. There are two cases to consider 
depending on whether $\widetilde{D}=\{x\}$ or not. In the case 
where $\widetilde{D}=\{x\}$, we have that the required top $R$-cycle is the singleton
$\{x\}$. We now pass to the case where $\widetilde{D} \neq\{x\}$. 
It follows that $\{x\} \subset D$. But then, we have that for at least one $x_{0} \in \widetilde{D}$ it holds that $x_{0}Rx$. 
Put
$$
A_{x}=\left\{y \in \widetilde{D} \mid(x, y) \in \overline{R}\right\}.
$$
We show that $x_{0} \in A_{x}$. We first show that $A_{x} \neq \emptyset$. Suppose to the contrary that 
$A_{x}=\emptyset$. Then, for each $y \in \widetilde{D}, (x, y) \notin \overline{R} \supseteq R$. It follows that 
$\widetilde{D} \backslash\{x\} \subset \widetilde{D}$ is an $R$-undominated subset of $X$, 
a contradiction because of the minimal 
character of $\widetilde{D}$. Let $\widetilde{D}(x)=\widetilde{D} \backslash A_{x}$. 
We now show that $\widetilde{D}(x)=\emptyset$. We proceed by contradiction, 
so let us assume that $\widetilde{D}(x) \neq \emptyset$. Then, for each $t \in A_{x}$ and each 
$s \in \widetilde{D}(x)$ we have 
$(t, s) \notin R$ for suppose otherwise, $(t, s) \in R$ implies that $(x, s) \in \overline{R}$ 
contradicting $s \in \widetilde{D}(x)$. 
Therefore, $\widetilde{D}(x) \subset \widetilde{D}$ is an $R$-undominated subset of $X$, 
which is again a contradiction. 
Hence, $\widetilde{D}(x)=\emptyset$ which implies that $A_{x}=\widetilde{D}$. Since $x_{0} \in \widetilde{D}$, 
we conclude that 
$(x, x_{_0}) \in \overline{R}$. Therefore, $x_{_0}, x$ belong to an $R$-cycle. 
Since $\widetilde{D}$ is a chain, for each $w, z\in \widetilde{D}$ we have that $(w,z)\in R$ or 
$(z,w)\in R$. Therefore, as in the case of $x_{_0}, x$, we conclude that $z, w$
belong to an $R$-cycle. It follows that $\widetilde{D}$ is a top $R$-cycle.
\par\smallskip\par\noindent
{\it Case {\rm (ii)}}: In this case we have $\widetilde{D}=\emptyset$. It follows that $D=\widetilde{G}$. That is $\widetilde{G}$ 
be the set of all $P(\overline{R})$-chains $G$ in $D$. Let $\preceq_{_{P(\overline{R})}}=P(\overline{R})\cup \Delta$ where
$\Delta=\{(x,x)\vert x\in X\}$.
Clearly, $\preceq_{_{P(\overline{R})}}$ is a partial order.
There holds that $x^{\ast}\in X$ is a $P(\overline{R})$-maximal element if and only if 
$x^{\ast}\in X$ is a $\preceq_{_{P(\overline{R})}}$-maximal element.
If
$\widetilde{G}=\{x\}$ for some
$x\in X$, then $x$ is $\preceq_{_{P(\overline{R})}}$-maximal element and thus $\{x\}$ is a top $R$-cycle.
Suppose that $X$ contains at least two elements. Let $\mathcal{U}$ be the family  
of all conforming sets in  $\widetilde{G}$ with respect to $\preceq_{_{P(\overline{R})}}$.
This family is not empty. Indeed, let $x,y$ be two distinct elements of $\widetilde{G}$ such that
$y\preceq_{_{P(\overline{R})}}x$ and $x\neq y$. Then, $\preceq_{_{P(\overline{R})}}$ is a well order of the set $\{y\}$ 
and for $y\in \{x\}$ we have $x=f(Q(\{y\}, x))$ where $f$ is a choice function which assigns to every chain 
$\mathfrak{C}$ of $\widetilde{G}$, a strict $\preceq_{_{P(\overline{R})}}$-upper bound $f(\mathfrak{C})$. According to the main result of \cite{lew} 
we have that
$\widetilde{G}$ has a $\preceq_{_{P(\overline{R})}}$-maximal element $x^{\ast}$. 
This is a 
$P(\overline{R})$-maximal element as well. 
If $x^{\ast}$ is an $R$-undominated element, then $\{x^{\ast}\}$ is a top $R$-cycle. 
 If this is not the case, then
there exists $\widetilde{y}\in D$ such that $(\widetilde{y},x^{\ast})\in R$.
Since $x^{\ast}$ is a $P(\overline{R})$-maximal element we
have that $(\widetilde{y},x^{\ast})\notin P(\overline{R})$. 
It follows that
$(x^{\ast},\widetilde{y})\in \overline{R}$ and $(\widetilde{y},x^{\ast})\in \overline{R}$.
Let
\begin{center}
$A_{{x^{\ast}}}=\{y \in D \vert (x^{\ast},y)\in \overline{R}$ and $(y,x^{\ast})\in \overline{R}\}$.
\end{center}
Since $\widetilde{y}\in A_{{x^{\ast}}}$ we have that $A_{{x^{\ast}}}\neq \emptyset$.
We show that $A_{{x^{\ast}}}$ is a top $R$-cycle.
Suppose to the contrary that 
there exists $t\in D\setminus A_{{x^{\ast}}}$ sush that $(t,s)\in R$ for some $s\in A_{{x^{\ast}}}$.
Then, $(s,x^{\ast})\in \overline{R}$ implies that $(t,x^{\ast})\in \overline{R}$. Since $x^{\ast}$ is a $P(\overline{R})$-maximal 
element we
have that $(t,x^{\ast})\notin P(\overline{R})$. Therefore, $(x^{\ast},t)\in \overline{R}$ which jointly to $(t,x^{\ast})\in \overline{R}$
implies that $t\in A_{{x^{\ast}}}$, a contradiction.
The last contradiction shows that $A_{{x^{\ast}}}$ is a top $R$-cycle which completes the proof.

\par\smallskip\par\noindent
(ii)$\Rightarrow$(i).\ 
Suppose that the generalization of Zorn lemma holds. Let 
$\mathcal{F}\subseteq \mathcal{P}(X)$
be a collection of nonempty sets. 
Let also $R^{\ast}$ be the quotient relation of $R$ in $X$. Then, $R^{\ast}$ is a partial order in 
$X_{_\approx}$.  Since $X$ has a top $R$-cycle we conclude that 
$X_{_\approx}$ has an $R^{\ast}$-maximal element. Therefore, according to the classical Zorn lemma, 
there exists a choice function for every collection of non-empty sets in
$X_{_\approx}$. 
Let $\mathcal{F}=\{F\vert F\subseteq X\}$ be a collection of nonempty sets in $X$ and let 
$[\mathcal{F}]=\{[F]\vert F\in \mathcal{F}\}$ be the corresponding collection in the quotient space 
$(X_{_\approx}, R^{\ast})$.
Suppose that 
$\widetilde{f}$ is a choice function for $[\mathcal{F}]$ in $X_{_\approx}$. Then, if $[F]\in [\mathcal{F}]$ 
we have $\widetilde{f}([F])\in [F]$. That is,
$\widetilde{f}([F])=[x]$ for some $[x]\in [F]$. 
Let  $x\in F$ be the representative of $[x]$. Then, from $\pi(F)=[F]$, $\widetilde{f}([F])=[x]$ and 
$\pi^{-1}([x])=x$ we conclude that $\pi^{-1}(\widetilde{f}(\pi(F)))=x\in F$. Therefore, 
$\pi^{-1}\circ \widetilde{f}\circ \pi$
is a choice function in
$X$.

$\qedsymbol$

\par\noindent
If the binary relation $R$ is a partial order, then it is evident that a top $R$-cycle is a singleton, that is, 
a maximal element. Thus, the above theorem generalizes Zorn's original result. 
\par\noindent
 As we mentioned at the beginning of this note, an example of the importance of Zorn's extended lemma 
which we introduce and can be applied in social choice theory is the theory of general solution concepts.
One of the most important general solution concepts is the Schwartz set which is defined in \cite{sch}.
This is done as follows: A choice function in social choice theory usually
represents those alternatives chosen by the individual or society. The traditional 
choice-theoretic approach takes behavior as rational if there is a binary relation 
$R$ such that for all non-empty subsets of $X, \mathfrak{C}(A)=\mathcal{M}(A, R)$ where $\mathcal{M}(A, R)$ is 
the set of $R$-maximal elements
of $A$. To deal with the case where the set of maximal elements is empty, Schwartz has proposed 
the general solution concept of $\mathcal{G}\mathcal{O}\mathcal{C}\mathcal{H} \mathcal{A}$ described above which is
equivalent to the union of all $R$-undominated elements and all top $R$-cycles in $X$.

\par\noindent
The following corollary of the above theorem provides an existence result for the Schwartz set which 
is equivalent to the extended Zorn's lemma. The proof of this corollary is obvious.
\par\smallskip\noindent
{\bf Corrolary}. Let $R$ be a binary relation on a set $X$. Assume that each chain
has an upper bound in $X$. Then, the Schwartz set is non-empty.
\par\smallskip\noindent
{\it Proof}. According to the above theorem the space $X$ has a minimal $R$-undominated set $D$.
Then, from \cite[Deb's theorem]{and} 
we have that $X$ has an $R$-maximal element or 
a strong top $R$-cycle.
It follows that the Schwartz set is non-empty.

\end{document}